\documentclass{article}
\usepackage{amssymb}
\usepackage{amsmath}
\usepackage{amsthm}

\numberwithin{equation}{section}

\begin{document}

\title{Constructing completely integrable fields by a generalized-streamlines method}
\author{Antonella Marini$^1$ and
Thomas H. Otway$^2$
\\ \\
\textit{$^1$Dipartimento di Matematica, Universit\`{a} di
L'Aquila,}\\ \textit{67100 L'Aquila, Italy } \\ \textit{$^{1,2}$Department of Mathematical Sciences,
Yeshiva University, }\\ \textit{New York, New York 10033}}
\date{}
\maketitle

\begin{abstract}

The classical approach to visualizing a flow, in terms of its streamlines, motivates a topological/soft-analytic argument for constrained variational equations. In its full generality, that argument provides an explicit formula for completely integrable solutions to a broad class of $n$-dimensional quasilinear exterior systems. In particular, it yields explicit solutions for extremal surfaces in Minkowski space and for Born--Infeld models.

\medskip

\end{abstract}

\section{Introduction}

The construction of streamlines associated to vector fields which arise as solutions to Laplace's equation goes back at least to Faraday and continues to be studied in every elementary physics laboratory, supplemented in more advanced courses by interpretation in the language of analytic functions and their conformal mappings; see, \emph{e.g.}, Ch.\ 1 of \cite{Me}. The corresponding literature for streamlines associated to stationary, irrotational, isentropic compressible flow is hardly less familiar. Notable, for example, is the phase-space treatment in Sec.\ 110 of \cite{CF}, the description of stream tubes in Sec.\ 7.1 of \cite{Ch}, and the treatment in Sec.\ 8 of \cite{Be} in the language of quasiconformal mappings.

In this paper we interpret the analytic aspects of streamlines in a different way. We focus on the role of the stream function, and of its topology, in inverting the differential operator for the continuity equation for compressible flow. The initial discussion proceeds via conventional vector analysis; but we will eventually find it useful (in Sec.\ \ref{Hdg}) to introduce invariant language. Just as Laplace's equation has a generalized interpretation in terms of Hodge theory, in which sections of the tangent bundle are replaced by sections of an exterior power of the cotangent bundle, nonlinear continuity equations have a generalized interpretation in terms of a \emph{nonlinear} Hodge theory \cite{SS1}. Thus geometric variational models may be expressed in terms of ``generalized streamlines" via a choice of exterior power for the cotangent bundle and a choice of the scalar function $\rho$ representing mass density in the original interpretation. This broader understanding converts the conventional stream function into a versatile tool for constructing explicit solutions in a variety of physical and geometric contexts. Note that in the past, Hodge-theoretic techniques seem to have been applied to continuity equations in the elliptic (subsonic) region only, whereas our methods extend to both sides of the sonic boundary. In particular, in Sec.\ \ref{secpatch} we apply our method of solution to a continuity equation in both the elliptic and hyperbolic regions individually, and paste the two solutions together along the sonic boundary.

\subsection{An outline of the method} \label{secintro}

Although there are many \emph{ad hoc} arguments that work for particular equations, few methods are available for producing explicit solutions to entire classes of quasilinear partial differential systems. This is especially the case if the system changes from elliptic to hyperbolic type along a smooth hypersurface in its domain. We begin by considering the exact solvability of systems having the form
\begin{equation} \label{NLH-intro}
\nabla\cdot\left(\rho\left(Q\right)\mathbf{w}\right)=0,
\end{equation}
\begin{equation} \label{exact}
\nabla\times\mathbf{w}=0,
\end{equation}
where $\mathbf{w}\equiv (w_1, \dots, w_n)$ is an unknown vector-valued function on a domain $\Omega\subset \mathbb R^n;$ $Q\equiv\vert\mathbf{w}\vert^2;$ $\rho\,:\, Dom_\rho\subset\mathbb{R}^+\cup\left\{0\right\}\rightarrow\mathbb{R}^+\cup\left\{0\right\}$ is a prescribed, continuously differentiable function on its domain of definition. By eq. \eqref{exact} we mean, by perhaps a slight abuse of notation, the  system of $n\times (n-1)$ equations
$$\partial_i w_j -\partial_j w_i = 0 \,,\;\mbox{ for all } i\neq j \in \{1, \dots n\}\,,$$
equivalent to the requirement that all the $2\times 2$ minors of the $2\times n$ matrix
\begin{equation}
\label{minors}
\left(
\begin{array}{cccc}
\partial_1 &\partial_2 &\dots &\partial_n\\
w_1 & w_2 &\dots & w_n
\end{array}
\right)\,,
\end{equation}
vanish.
Equations (\ref{NLH-intro}, \ref{exact}) are a vectorial special case (corresponding to the case of differential 1-forms in $n$ dimensions) of the \emph{nonlinear Hodge equations} introduced in \cite{SS1}. The system is satisfied by a large number of physical and geometric models, including the irrotational steady flow of a compressible fluid; light near a caustic; shallow water hydrodynamics; ball lightning in the nonlinear conductivity model; and fields which can be represented as a non-parametric extremal surface in $\mathbb{R}^3$ or $\mathbb{M}^3.$ See Sec.\ 2.7 and Chs.\ 5 and 6 of \cite{O4}, as well as references cited therein, for discussions of these applications.

Equation (\ref{NLH-intro}) is a variational equation associated to the energy functional
\[
E = \frac{1}{2}\int_\Omega\int_0^Q\rho(s)ds*1\,,
\]
while equation (\ref{exact}) constrains the resulting field to be locally conservative. In the language of elementary differential equations, eq.\ (\ref{exact}) implies an \emph{exactness condition} on contractible domains: if (\ref{exact}) is satisfied, then we can write
\begin{equation} \label{cons}
\mathbf{w} = \nabla \zeta
\end{equation}
for some scalar function $\zeta.$

We develop a method by which we derive the solution formula for the variational equation (\ref{NLH-intro}). This formula holds globally on contractible domains and locally on more general domains. A less specialized version of the method, in which streamlines appear only locally, yields a global solution formula on non-contractible domains as well. Extensions to boundary value problems by means of De Rham cohomology and Hodge decomposition theorems on manifolds are pursued in \cite{MO2}, Sec.\ 5.

\smallskip

Naturally,  condition (\ref{exact}) is not generally satisfied. So, continuing to follow the analogy of elementary differential equations, we look for an integrating factor. That is, we replace the constraint  (\ref{exact}) by the \emph{Frobenius condition}
\begin{equation} \label{frob-2}
\nabla\times\mathbf{w}=\mathbf{G}\times\mathbf{w},
\end{equation}
where $\mathbf{G}=\mathbf{G}(\mathbf x)$ is a vector function.  Again borrowing (and extending) the standard notation for the cross product of vector fields in $\mathbb{R}^3$, by eq. \eqref{frob-2} we mean the system of $n\times (n-1)$ equations
$$\partial_i w_j-\partial_j w_i = G_i w_j - G_j w_i\,,\; \mbox{ for all } i\neq j\in\{1,2,\dots, n\}\,,$$ obtained by equating the $2\times 2$ minors of the $2\times n$ matrix \eqref{minors} to the corresponding minors of
$$
\left(
\begin{array}{cccc}
G_1 & G_2 &\dots &G_n\\
w_1 & w_2 &\dots & w_n
\end{array}
\right).
$$
As a consequence of the Frobenius theorem \cite{E}, a solution to (\ref{frob-2}) in the context just described can always be written locally in the form
\begin{equation}
\label{intfact}
e^{-\eta}\mathbf{w} = \nabla\zeta,
\end{equation}
where $\eta$ and $\zeta$ are functions and $\mathbf{G}$ can always be chosen to be conservative:
\begin{equation}
\label{exactvectorial}
\mathbf{G} = \nabla\eta\,.
\end{equation}
Thus one of the goals of the method is to construct explicit solutions to eq.\ (\ref{NLH-intro}) which are \emph{completely integrable fields}. By this we mean that they will possess an \emph{integrating factor}, namely a function $\eta$ satisfying \eqref{intfact} or equivalently, in this vectorial case, satisfying (\ref{frob-2}), \eqref{exactvectorial}.
See \cite{MO1}, Secs.\ 2.1, 2.2 for a discussion.

Notice that the aim of our method is not to solve (\ref{NLH-intro}, \ref{frob-2}) for  \emph{prescribed} $\mathbf{G},$ including the original case $\mathbf{G}=0.$ Indeed, solving (\ref{NLH-intro}, \ref{frob-2}) for prescribed $\mathbf{G}$ would entail imposing the condition that our explicit solution also solve an additional non-linear differential equation. However, a vector field $\mathbf{G}$ for which $\mathbf w$ satisfies \eqref{frob-2} can be computed once the explicit solution to eq.\ (\ref{NLH-intro}) has been found. Moreover, as an artifact of the proof, we derive the auxiliary equation that would need to be satisfied in order to solve the system for prescribed $\mathbf{G}.$

Analogous arguments apply to the case in which one is considering eqs. \eqref{vareq-k=n-1}, \eqref{frob-k=n-1}, below. The solutions to system \eqref{vareq-k=n-1} found by our formula  are also completely integrable, as they ordinarily satisfy the Frobenius condition \eqref{frob-k=n-1} for some vector function \eqref{exactvectorial}.

\smallskip

We apply the method to various typologies of systems of equations, where the typology is determined by  two different dimensional parameters, $k$ (representing the degree of the associated differential form) and $n$ (representing the dimension of the domain).
In all cases, the method decomposes into two steps: a first step yielding  the normalized solutions,  and a second step yielding the magnitude of the vector solutions $Q(\mathbf x)$.
A feature of the method is that the two steps are essentially unrelated. In particular, the first step is independent of the density function $\rho$  and depends only on the dimensional parameters $k$ and $n$; the second step is totally independent of the system under consideration and
is concerned almost exclusively with the density function $\rho$, via the invertibility of the real-valued function of one variable $\xi = \phi(t)\equiv t\rho^2(t),$ bearing in mind that the analysis of the singular sets will depend on the dimension of the domain.

\smallskip

The system (\ref{NLH-intro}, \ref{frob-2}) was introduced for the general case of differential forms of arbitrary order  in \cite{O1} and studied in \cite{O2}; it received an extensive treatment in \cite{MO1}. We call that generalized form of (\ref{NLH-intro}, \ref{frob-2}) the \emph{nonlinear Hodge--Frobenius equations}. Equation \eqref{frob-2}  represents the mildest natural weakening of the
side  condition on the variations in the conventional nonlinear Hodge equations; compare also with the conservative system studied in \cite{ISS}.
The Frobenius condition (\ref{frob-2}) or \eqref{frob-k=n-1}, and the special cases (\ref{exact}) or \eqref{exact-k=n-1} respectively, corresponding to the choice $\textbf{G}=0$, generally have physical or geometric significance. For example, if the mass density $\rho$ in eq.\ (\ref{NLH-intro}) satisfies
\begin{equation} \label{shallow}
\rho(Q) = 1-\frac{Q}{2}\,,
\end{equation}
then we obtain  a model for steady, shallow, hydrodynamic flow under an appropriate normalization; see, \emph{e.g.}, Sec.\ 10.12 of \cite{St}. If  (\ref{exact}) is satisfied, then the flow is irrotational; if (\ref{frob-2}) is satisfied, then the flow is vorticial (Secs.\ \ref{SS-shallow} and \ref{SS-special}). If we choose in \eqref{NLH-intro}  the density
\begin{equation}\label{extr}
\rho(Q)=\vert 1-Q\vert^{-1/2}, \,\, Q\ne 1\,,
\end{equation}
then we obtain a model for extremal surfaces, in 3-dimensional Minkowski space $\mathbb{M}^3,$ which can be expressed locally as the graph of a function $\zeta;$ see Sec.\ 1 of \cite{SSY}, Ch.\ 6 of \cite{O4}, and the references cited therein. Replacing (\ref{exact}) by (\ref{frob-2}) dilates the gradient vector of $\zeta$ without changing the direction of that vector (Sec.\ \ref{SS-extremal}). For other examples, see, Secs.\ 1--4 of \cite{MO1}.

\smallskip
\subsection{Organization of the paper}
In Sec.\ \ref{S-stream2} we develop the method in $2$ dimensions for the system  (\ref{NLH-intro}, \ref{frob-2}).

In Sec.\ \ref{S-streamn} we extend the method to $n$ dimensions, introducing the generalized stream flow generated by
the stream matrix $F$.

In Sec.\ \ref{S-k=n-1} we apply the method to a related system of $n(n-1)/2$ equations in $n$ dimensions.

In Sec.\ \ref{S-k} we give the general formulation of the method for differential forms of arbitrary degree on domains of arbitrary dimension, and apply the method to the Born--Infeld model.

In Sec.\ \ref{S-exa} we present examples and applications. These include an illustration, in Sec.\ \ref{secpatch}, of the important issue of how to patch the individual solutions in the hyperbolic and elliptic regions into a single solution which crosses the sonic transition with an acceptable degree of regularity.

\section{The method in two dimensions}
\label{S-stream2}
Denote by $\Omega$ a simply connected domain in $\mathbb R^2$. Applying the Poincar\'e Lemma to eq.\ (\ref{NLH-intro}) indicates that there is a scalar differentiable function $f\left(x,y\right)$ such that
\begin{equation} \label{poinc}
\rho(Q)\,\mathbf{w} =  \nabla_\perp f\,,
\end{equation}
having denoted by $ \nabla_\perp f$ the \emph{transverse gradient} of $f$ given by
\[
\nabla_\perp f \equiv -\left(\partial_y f\right) \hat{\imath}+ \left(\partial_x f\right) \hat{\jmath}.
\]
In the language of fluid dynamics, eq.\ (\ref{poinc}) asserts that there exists a \emph{stream function} $f$ for the divergence-free vector field $\rho(Q)\mathbf{w}.$ Conversely, for any prescribed scalar function  $f\left(x,y\right)$ defined on a  domain $\tilde\Omega$, the field $\nabla_\perp f(x,y)$ is divergence-free. So \eqref{poinc} leads to solutions of \eqref{NLH-intro}, even if $\tilde\Omega$ is not simply connected, in which case one may
prescribe $f$ with point singularities; \emph {cf.} also Remark $iv)$ at the end of this section. In seeking explicit solutions $\mathbf{w}$ of eq.\eqref{NLH-intro} via \eqref{poinc}, one encounters two problems. The first is that eq. \eqref{poinc} may not be easily inverted to yield $\mathbf{w}$ in terms of $f$; the second is that  $\rho$ may vanish for values of $Q$ attained on a subset $\gamma_0\subset\Omega.$ Both problems arise in actual applications. In what follows we illustrate how to proceed in order to invert eq. \eqref{poinc} and argue that in many cases the second problem does not occur, or is only apparent in the sense that it corresponds to the existence of removable singularities.

Our initial observation is that, because $\rho$ is a scalar function, eq.\ \eqref{poinc} yields for $\rho\ne0$ the precise expression in terms of $f$ for the normalized  solution
 \begin{equation}\label{normw}
  \hat{\mathbf{w}} \equiv \frac{\mathbf{w}}{|\mathbf{w}|}= \pm \frac{\nabla_\perp f}
  {\vert\nabla f\vert}.
  \end{equation}
The mass density $\rho$ does not appear explicitly in eq. \eqref{normw}.  Nevertheless, this expression holds only in the limit (or not at all) on a singular set on which $\nabla f = \mathbf 0.$ Such a set would itself depend on $\rho$ and $f,$ as will be illustrated later in the present section and in the examples of Sec.\ \ref{S-exa}.
An immediate consequence of our observation is that, in order to find a complete expression for $\mathbf{w},$ it suffices to find an equation satisfied by $\vert\mathbf{w}\vert$ alone in terms of $f$. To this end we take the squared norm on both sides of eq. \eqref{poinc} and find that the real-valued function of a real variable
\begin{equation} \label{phi}
\phi (t) \equiv t\rho^2(t),\,\,t\geq 0\,,
\end{equation}
can be usefully inverted to yield $t=Q\equiv \vert\mathbf{w}\vert^2.$

We say that the function $\phi$ restricted to an interval $I_1$   or  $I_2$ is  \emph{invertible with inverse of type} $1$,  or \emph{invertible with inverse of type} $2$, respectively, if
\[
\frac{d\phi}{dt}(Q)>0 \,,\; Q\in I_1\,, \quad \mbox {or } \;\frac{d\phi}{dt}(Q)<0\,,\; Q\in I_2\,.
\]
We denote by $\psi_1$,  $\psi_2$ the corresponding inverses.
The first choice corresponds to a region of $\Omega$ on which the system (\ref{NLH-intro}, \ref{frob-2}) is elliptic; the second corresponds to a region of $\Omega$ on which the system (\ref{NLH-intro}, \ref{frob-2}) is hyperbolic, as  can be easily computed.  In various models there are non-empty intervals  $I_1$ and $I_2$ on which $\phi$ is invertible with individual inverses in the elliptic and hyperbolic regimes, although we may not know \emph{a priori} the corresponding regions in the domain $\Omega$. Sometimes elliptic and hyperbolic solutions can be patched together along the  \emph{sonic curve} dividing the elliptic from the hyperbolic regime. We will illustrate some of these models in Sec.\ \ref{S-exa}.

Because $\rho$ is given and we have, by (\ref{poinc}),
\begin{equation}\label{gradf}
\phi(Q) \equiv Q\rho^2(Q)=\vert\nabla f\vert^2,
\end{equation}
 the quantity $Q$ is expressible in terms of the prescribed \emph{generalized streamline} $f$ (more precisely, in terms of the magnitude of its
gradient) via the inverse, or inverses, of $\phi.$ So \emph{a posteriori} the
elliptic and hyperbolic regimes  correspond to regions
of the plane where one uses inverses of type $1$ or of type $2$, respectively.   Using  \eqref{poinc}, \eqref{gradf} to solve for $\mathbf{w}$ in terms
of $f,$ we obtain solutions to \eqref{NLH-intro} expressed by
\begin{equation} \label{vareqsol}
\mathbf{w} = \frac{\nabla_\perp f}{\rho\left(\psi\left(\vert\nabla f\vert^2\right)\right)}\,,
\end{equation}
where $\psi$ may denote inverses on different monotonicity intervals.
These can be
defined on a subdomain
\[\Omega_f \equiv \{ (x,y)\in\Omega\cap Dom_{\nabla f} \;:\;
\vert\nabla f\vert^2\left(x,y\right)\in \mbox{Im}\,\phi \}\subseteq\Omega\,,
\]
except possibly on a singular set $\mathcal S\subset \gamma_s;$ here
$$\gamma_s\equiv\{(x,y)\in\Omega_f\;:\;
\phi^\prime\left(\psi\left(\vert\nabla f\vert^2\left(x,y\right)\right)\right)=0 \mbox { or  is undefined}\}\,,$$
where one may switch between different inverse functions of $\phi.$
We observe that the set $\gamma_s$ always contains the set
$$\gamma_0\equiv\{(x,y)\in\Omega_f\;:\;
\rho\left(\psi\left(\vert\nabla f\vert^2(x,y)\right)\right)=0 \mbox{ and }\psi\left(\vert\nabla f\vert^2(x,y)\right)\neq 0\}\,,$$
on which the solutions \eqref{vareqsol} may blow up or be undetermined. That is, $$\gamma_0 \subset\gamma_s\,,$$
and such an inclusion may be proper.
This follows easily from the relation $$\phi^\prime(Q) = \left(Q\rho^2(Q)\right)^\prime=\rho(Q) \left(\rho(Q) +2Q\rho^\prime(Q)\right).$$
Moreover,
$\gamma_0\subset
\{(x,y)\in\Omega_f\;:\;\vert\nabla f\vert^2\left(x,y\right)=0\}$.
In fact, by \eqref{gradf} -- and its alternate form \eqref{altern0}, see below -- the condition $\rho\left(\psi\left(\vert\nabla f\vert^2(x,y)\right)\right)=0$ implies $\vert\nabla f\vert^2\left(x,y\right)=0$. This inclusion is proper when there are points $(x,y)\in \Omega_f$ at which $\psi\left(\vert\nabla f\vert^2(x,y)\right)=0$ and  $\rho\left(\psi\left(\vert\nabla f\vert^2(x,y)\right)\right)\neq 0.$
Also note that, by the alternate form of \eqref{vareqsol} given by \eqref{altern} below, it is clear that $\mathbf w$ is defined and equals zero at points where $\psi(\vert\nabla f\vert^2(x,y))=0,$ even when $\rho\left(\psi\left(\vert\nabla f\vert^2(x,y)\right)\right)=0$ at these points.

Although for a smooth prescribed function $f$,  it may be possible to define a corresponding $\mathbf w$ in such a way to transition with continuity between different inverses $\psi$ -- in particular, between elliptic and hyperbolic regimes -- it will in general not be possible to do so with higher regularity. So,  although $\mathbf w$ may be defined with some regularity on the set $\Omega_f - \gamma_0$,  eq. \eqref{NLH-intro} may not hold in the classical sense across the set $\gamma_s$. For this reason, we regard $\gamma_s$  as a possibly singular set.
In important applications, this set is a smooth curve -- the sonic curve mentioned earlier.

For analogous reasons, one also regards the set $\gamma_\infty$, defined as
$$\gamma_\infty \equiv\{(x,y)\in \Omega \,:\, \rho\left(\psi\left(\vert\nabla f\vert^2(x,y)\right)\right) \mbox { is undefined }\},
$$ as a singular set for the equation. As $f$ is prescribed, one can avoid altogether the sets $\gamma_0$ and $\gamma_s$,  and produce examples for which $\gamma_0=\gamma_s=\emptyset$  and ${\mathbf w}$ is smooth. Nevertheless it is more interesting, mathematically and in terms of the applications, to produce examples for which $\gamma_s$ and perhaps also $\gamma_0$  are non-empty.

Strictly speaking,  the set of points on which $\psi(\vert\nabla f\vert^2)$ becomes unbounded does not belong to $\Omega_f$, and the definition of ${\mathbf w}$ cannot be extended to these points, as is apparent from the alternate expression for ${\mathbf w}$ given by \eqref{altern}, below.

For a given function $f$, one has $\Omega_f = \Omega$ -- that is, the corresponding solutions $\mathbf{w}$  live on all of $\Omega$, except possibly on a singular  set $S\subset \gamma_s$ -- if and only if
\[
\Sigma_f\equiv \left\{\vert\nabla f\vert^2\left(x,y\right):\left(x,y\right)\in\Omega\right\}\subset \mbox{Im} \,\phi.
\]
In the particular case in which
\[
\Sigma_f\subsetneq \mbox{Im}\,\phi\,,
\]
the full range of $\phi$ is not available and one may be able to invert only in the elliptic regime, or only in the hyperbolic regime, for that particular choice of $f.$
There are also choices of $f$ for which $\Omega_f=\emptyset.$
The example of the density function  $\rho= c/\sqrt Q,$ where $c$ is a constant, illustrates this possibility when $\vert\nabla f\vert \neq c$.
In fact, for this density function one has
$$\frac{c}{\sqrt Q} \,\mathbf w = \nabla_\perp f\,,$$
which admits solutions if and only if $\vert\nabla f\vert(x,y) = c$ $\,\forall \, (x,y)\in \Omega.$ Notice that only the normalized $\mathbf w /\sqrt Q$ is specified by this equation, while the function $Q(x,y)$ can never be determined. In fact in this example  $\phi(Q)\equiv c^2$ cannot be inverted and one always has $\emptyset=\gamma_0 \subset \gamma_s=\Omega.$
In the case $\vert\nabla f\vert (x,y) = c$ $\,\forall \,(x,y)\in \Omega,$ one has $\Omega_f =\Omega$ and the normalized solution corresponding to $f$ is defined everywhere. Multiplying by the arbitrary function $\sqrt{Q(x,y)},$ one obtains a family of solutions. This dramatically illustrates the importance of the range of $\vert\nabla f\vert^2$ in determining solutions to eq.  \eqref{NLH-intro}.

The simple example of the density function $\rho=Q^{-\frac{1}{4}}$ illustrates instead a case for which $\gamma_\infty \neq \emptyset$, and yet one can still carry out our method and obtain solutions  $$\mathbf w = \nabla_\perp f \;\vert \nabla f \vert\,,$$
defined on the domain of $\nabla f$ and smooth for smooth functions  $f$. In this example, $\phi = \sqrt Q$ is everywhere defined and increasing. Nonetheless, although the solutions defined above may be everywhere smooth, they are technically solutions to
\eqref{NLH-intro} only on  $\{ (x,y)\in\Omega\,:\, \nabla f (x,y)\neq \mathbf 0\}.$

\bigskip
\textbf{Remarks}.

\medskip
\emph{i}) In (\ref{vareqsol}) we have divided by $\rho$ -- which may vanish -- but, as mentioned
previously, this problem is only apparent in many cases. In fact,
suppose that $\rho(Q) =0$ for some finite value(s) of $Q$, say for
example $\rho(Q_1) = 0$. This gives \emph{a posteriori} that
$\vert\mathbf{w}\vert = \sqrt{Q_1}$ is bounded at the point(s) $p_1\in\Omega$ at which the value $Q_1$ is attained, despite the vanishing of $\rho$ at $p_1.$ It is possible, however, that in some cases ${\mathbf w}$ may be not well defined at $p_1$; \emph{cf}.\ Remark \emph{ii}, below. That is the case
for example in some models in which $\mathbf{w}$ develops point
singularities. It is also possible that $\rho$ will vanish at points
of $\Omega$ for which $Q$ becomes infinite. That is the case, for
example, in the variational model considered in Sec.\ \ref{SS-extremal}, when one
chooses to invert in the hyperbolic regime; see eq.\ \eqref{wt}
for $\vert\nabla f\vert\to 1.$

\medskip
\emph{ii}) Rewriting \eqref{gradf} as
\begin{equation}
\label{altern0}
\vert\nabla f\vert^2= \phi(\psi(\vert\nabla f\vert^2))=\rho^2(\psi(\vert \nabla f\vert^2))\,\psi(\vert \nabla f\vert^2)\,,
\end{equation}
we obtain an alternate expression for $\mathbf{w}$ which may be used in place of \eqref{vareqsol} for finding  explicit solutions to \eqref{NLH-intro}, namely,
\begin{equation}
\label{altern}
\mathbf{w} = \pm \frac{\nabla_\perp f}{\vert\nabla f\vert}\,\sqrt{\psi\left(\vert\nabla f\vert^2\right)}\,.
\end{equation}
This illustrates what may happen at points of $\gamma_0:$  even though $Q$ may remain bounded at these points, in various examples and applications the normalized vector function ${\mathbf w}/|{\mathbf w}|$ may not be defined (even as a limit) at these points.

Formally, the two formulas for $\mathbf{w}$, \eqref{vareqsol} and \eqref{altern},  coincide if the minus sign in (\ref{altern}) is assigned to regions where $\rho\left(\psi(\vert \nabla f\vert^2)\right)$ is negative (generally excluded in the applications). On the other hand,  the minus sign in (\ref{altern}) may be avoided altogether, as
\[
\vert\nabla f\vert\left(x,y\right)=\vert\nabla \left(-f\right)\vert\left(x,y\right)\,;
\]
the choices of either $-f$ or $f$ in \eqref{altern} are equally legitimate and yield the same range: $\Sigma_{-f}=\Sigma_f.$

\medskip
\emph{iii}) The vector function $\mathbf{w}$ given by eq.\ (\ref{vareqsol}) or \eqref{altern} satisfies (\ref{NLH-intro}) but may not satisfy (\ref{exact}). However, for any given $\mathbf{G},$ including the case $\mathbf{G}=0,$ eq.\ (\ref{frob-2}) will be satisfied with the choice
\begin{equation} \label{G-2d}
\mathbf{G} = -\mathbf{G}_1-\frac{\nabla\rho\left(\psi\left(\vert \nabla f\vert^2\right)\right)}{\rho\left(\psi\left(\vert\nabla f\vert^2\right)\right)},
\end{equation}
whenever $f$ satisfies the equation
\begin{equation} \label{G_1eq-2d}
\Delta f +\mathbf{G}_1\cdot\nabla  f=0
\end{equation}
for some continuously differentiable vector function $\mathbf{G}_1.$
Solving (\ref{G-2d}) for $\mathbf{G}_1$ and substituting the result into (\ref{G_1eq-2d}), we find that the latter is a nonlinear equation for $f$. In any case, our aim is not to solve eq.\ (\ref{G_1eq-2d}) for $f$ in terms of a prescribed $\mathbf G_1$, but to obtain completely integrable solutions.

Indeed, for the reasons explained in Sec.\ 1, we require only that condition (\ref{frob-2}) be satisfied for \emph{some} $\mathbf{G}$ once $\mathbf{w}$ has been determined by formula (\ref{vareqsol}). This can always be achieved, except on the (possibly empty) set
\begin{equation}
\label{nonG-2d}
\gamma_G\equiv\{ (x,y)\in \Omega_f\subset\mathbb{R}^2\,:\, \nabla f (x,y) =\mathbf 0\,; \;\Delta f(x,y) \neq 0\,\}\,,
\end{equation}
by choosing for example,
\begin{equation}\label{G_1ex-2d}
\mathbf{G}_1 = -\Delta f\left(\frac{\nabla f}{\vert\nabla f\vert^2}\right)\,.
\end{equation}
Because adding to $\mathbf G_1$ a vector function $\mathbf H$  satisfying $\mathbf H\cdot \nabla f =0$  has no effect on eq. \eqref{G-2d}, there will be, in general,  infinitely many $\mathbf G$ satisfying eq. \eqref{frob-2}. We know from the theory \cite{E} that one can specify $\mathbf G$ to be  a conservative vector field.

\medskip

\emph{iv}) The solution formula to \eqref{NLH-intro} found by the method of generalized streamline can be extended to domains $\tilde\Omega$ which are not simply connected, by replacing the transverse gradients $\nabla_\perp f\,$ in \eqref{poinc} by arbitrary divergence-free vector fields $\mathbf\alpha.$ Once this substitution is made, the remainder of the procedure does not change,  yielding
$$\mathbf w= \frac{\mathbf\alpha}{\rho(\psi(\vert\mathbf\alpha\vert^2))}\,.$$
A calculation analogous to the one done in Remark $iii)$ shows that the Frobenius condition \eqref{frob-2} for $\mathbf w$
 holds with $\mathbf G$ satisfying \eqref{G-2d} and $\mathbf G_1$ given by
$$
\mathbf{G}_1 = - (\partial_2\alpha_1-\partial_1\alpha_2 )\,\left(\frac{\mathbf\alpha}{\vert\mathbf\alpha\vert^2}\right)=- \nabla_\perp \cdot \mathbf\alpha\,\left(\frac{\mathbf\alpha}{\vert\mathbf\alpha\vert^2}\right)\,,$$
with $\nabla_\perp \equiv (-\partial_2, \partial_1).$

\section{The method in $n$ dimensions}
\label{S-streamn}
The method described in Sec.\ \ref{S-stream2} for the $2$-dimensional case can be extended to vectors in $n$ dimensions as follows.

The Poincar\'e Lemma applied to eq.\ (\ref{NLH-intro}) on contractible domains $\Omega\subset\mathbb R^n$ gives the existence of an $n\times n$ real matrix-valued function $F\equiv (f_{ij}),$ with $F$ skew symmetric ($f_{ii}=0$ and $f_{ij} = - f_{ji}$, $\forall i,j\in\{1, \dots, n\}$),  such that
\begin{equation} \label{poincn}
\rho(Q)\,\mathbf{w} =  \nabla_\perp F\,,
\end{equation}
where by analogy with Sec.\ \ref{S-stream2} we define the \emph{transverse gradient of a skew symmetric matrix $F$}, denoted  $\nabla_\perp F$, to be the $n$-vector
\[
\nabla_\perp F \equiv \left(\sum_{j\neq 1} \partial_j f_{1j} \right) {\hat{\imath}}_1+ \dots +  \left(\sum_{j\neq n} \partial_j f_{nj} \right) {\hat{\imath}}_n\,.
\]
(Notice that for $n=2$ this corresponds to the description in Sec.\ \ref{S-stream2} after making the substitution $f = -f_{12} = f_{21}.$)

Extending the language of fluid dynamics, eq.\ (\ref{poincn}) asserts that there exists a \emph{stream matrix} $F$ for the divergence-free vector field $\rho(Q)\mathbf{w}.$ On domains $\tilde\Omega$ which are not simply connected, eq.\ \eqref{poincn} still leads to solutions of \eqref{NLH-intro}. More generally, one can replace $\nabla_\perp F$ by  arbitrary sufficiently smooth divergence-free $n$-dimensional vector fields $\mathbf\alpha$ on the right-hand side of \eqref{poincn}.

In order to find explicit solutions $\mathbf{w}$ to eq. \eqref{NLH-intro} in terms of given matrices $F$, we continue to proceed by analogy with Sec.\ \ref{S-stream2}.

It is again possible to express the normalized solution to \eqref{NLH-intro} in terms of the stream matrix by the formula
 \begin{equation}\label{normwn}
  \hat{\mathbf{w}} \equiv \frac{\mathbf{w}}{|\mathbf{w}|}= \pm \frac{\nabla_\perp F}
  {\vert\nabla_\perp F\vert}\,,
  \end{equation}
in which the particular choice of mass density $\rho$ does not appear, although  the singular set for this formula may depend on $\rho$.

Because $\rho$ is given, eq.  (\ref{poincn}) yields as in the two-dimensional case an explicit relation between the two scalar quantities $|\nabla_\perp F|^2$ and $Q\equiv\vert {\mathbf w}\vert^2$, namely
\begin{equation}\label{gradfn}
\phi(Q) \equiv Q\rho^2(Q)=\vert\nabla_\perp F\vert^2\,,
\end{equation}
enabling us to express the quantity $Q$ in terms of the prescribed
stream matrix $F$ via the inverse(s) of $\phi.$

As in the two-dimensional case, the
elliptic and hyperbolic regimes for the system of eqs.\ \eqref{NLH-intro}, \eqref{frob-2} correspond to the regions
of $\Omega\subset\mathbb{R}^n$ on which one uses inverses of type $1$ or of type $2$.
Denoting simply by $\psi$ the various inverses of $\phi$ restricted to its monotonicity intervals,
one obtains
solutions to \eqref{NLH-intro} of the form
\begin{equation} \label{vareqsoln}
\mathbf{w} = \frac{\nabla_\perp F}{\rho\left(\psi\left(\vert\nabla_\perp F \vert^2\right)\right)}\,.
\end{equation}
These are defined on subdomains
\[\Omega_F \equiv \{\mathbf x\equiv(x^1, \dots, x^n) \in\Omega\cap Dom_{\nabla_\perp F} \;:\;
\vert\nabla_\perp F\vert^2 (\mathbf x )\in \mbox{Im}\,\phi \}\subseteq\Omega\,,
\]
except possibly on a set $\mathcal S\subset \gamma_s$ with
$$\gamma_s\equiv\{\mathbf x\in\Omega_f\;:\;
\phi^\prime\left(\psi\left( \vert\nabla_\perp F\vert^2 (x,y)\right)\right)= 0 \mbox { or  is undefined }\}\,.$$
We observe once more that the set
$$\gamma_0\equiv\{ \mathbf x\in\Omega_F\;:\;
\rho(\psi(\vert\nabla_\perp F\vert^2 (\mathbf x)))= 0\mbox { and } \psi(\vert\nabla_\perp F\vert^2 (\mathbf x))\neq 0\}\,,$$
where $\mathbf w$ may blow up or be undefined, is contained in $\gamma_s$.
The latter is to be regarded as a possibly singular set. Even though,  for a smooth prescribed stream matrix $F$ it may still be possible to define ${\mathbf w}$ as a vector function with continuity on the set dividing the elliptic from the hyperbolic regime,  eq. \eqref{NLH-intro} will not hold on this set.
In this context, under suitable regularity assumptions, the sets $\gamma_0$ and $\gamma_s$ will be smooth $(n-1)$-dimensional hypersurfaces.
As $F$ is prescribed, one can keep away from these singular sets and produce various examples for which $\gamma_0=\gamma_s=\emptyset$  and ${\mathbf w}$ is smooth. Nevertheless, it is again more interesting, intrinsically as well as in terms of possible applications, to produce examples for which $\gamma_s$, or  $\gamma_0$ and $\gamma_s,$ are non-empty.
Again, the set of points
for which $\psi(\vert\nabla F\vert^2)(x,y)$ becomes unbounded is excluded from $\Omega_f$ by definition, and ${\mathbf w}$ cannot be extended to these points, as is clear from the alternate expression for ${\mathbf w}$ given by
\begin{equation}
\label{alternn}
\mathbf{w} = \pm \frac{\nabla_\perp F}{\vert\nabla_\perp F\vert}\,\sqrt{\psi\left(\vert\nabla_\perp F \vert^2\right)}\,,
\end{equation}
an $n$-dimensional analogue of \eqref{altern}.

For a given stream matrix $F$, one has $\Omega_F = \Omega$ -- that is, the solutions $\mathbf{w}$ given above live on all of $\Omega$, except possibly on the sets $\gamma_0,$ on which  ${\mathbf w}$ may be undefined, and $\gamma_s$, on which ${\mathbf w}$ is not technically a solution to eq. \eqref{NLH-intro}  -- if and only if
\[
\Sigma_F\equiv \left\{\vert\nabla_\perp F \vert^2\left(x,y\right):\left(x,y\right)\in\Omega\right\}\subset \mbox{Im} \,\phi\,.
\]
In various models there are non-empty intervals $I_1$, $I_2$ on which $\phi$ is invertible with individual inverses on the elliptic and  hyperbolic regimes. Some of these models will be reviewed in Sec.\ 5.

\medskip
The analysis of the Frobenius condition is similar to the analysis carried out in Sec.\ \ref{S-stream2}. The Frobenius condition is obtained with
 \begin{equation} \label{G-nd}
\mathbf{G} = -\mathbf{G}_1-\frac{\nabla\rho\left(\psi\left(\vert \nabla_\perp F \vert^2\right)\right)}{\rho\left(\psi\left(\vert\nabla_\perp F \vert^2\right)\right)}\,,
\end{equation}
where  $\mathbf G_1$  satisfies  the system of $n(n-1)/2$ equations
\begin{equation} \label{G_1eq-nd}
\nabla \times \nabla_\perp F +\mathbf{G}_1\times \nabla_\perp F=0\,.
\end{equation}
On domains $\tilde\Omega$ which are not simply connected, in the more general context the Frobenius condition holds with  $\nabla_\perp F$ in eqs. \eqref{G-nd}, \eqref{G_1eq-nd} replaced by arbitrary sufficiently smooth divergence-free $n$-dimensional vector fields $\mathbf\alpha$.

\section{Extension of the method to systems}
\label{S-k=n-1}
The method also applies to the problem  in which eq. \eqref{NLH-intro} is replaced by the system of equations
\begin{equation}
\label{vareq-k=n-1}
\partial_i (\rho w_j) -\partial_j (\rho w_i )= 0 \,,\;\mbox{ for all } i\neq j \in \{1, \dots n\}\,,
\end{equation}
for a vector function ${\mathbf w}({\mathbf x}) = (w_1, \dots, w_n)({\mathbf x})$  on a domain $\Omega\subset\mathbb{R}^n$, for prescribed  $\rho=\rho (\vert {\mathbf w}\vert^2)$.
Equation \eqref{vareq-k=n-1} prescribes the vanishing of the  $2\times 2$ minors of the matrix
$$
\label{minors-rho}
\left(
\begin{array}{cccc}
\partial_1 &\partial_2 &\dots &\partial_n\\
\rho w_1 & \rho w_2 &\dots & \rho w_n
\end{array}
\right)$$
In this case, eq. \eqref{exact} and its more general version \eqref{frob-2} are replaced by
\begin{equation}
\label{exact-k=n-1}
\nabla\cdot {\mathbf w} \equiv \sum_{j=1}^n \partial_j w_j= 0
\end{equation}
and
\begin{equation}
\label{frob-k=n-1}
\nabla\cdot {\mathbf w} \equiv \sum_{j=1}^n \partial_j w_j= {\mathbf G}\cdot {\mathbf w} \equiv \sum_{j=1}^n G_j w_j
\end{equation}
respectively, where ${\mathbf G}({\mathbf x})\equiv (G_1, \dots, G_n)({\mathbf x})$ is a vector function.

Application of the method to this problem yields the solution formula for classical solutions on contractible domains \begin{equation} \label{sol-k=n-1}
{\mathbf w} =\frac{\nabla f}{\rho (\psi (\vert \nabla f\vert^2))}\,,
\end{equation}
or the alternate formula
\begin{equation}
\label{altern-k=n-1}
{\mathbf w} =\frac{ \nabla f}{\vert \nabla f\vert} \, \sqrt {(\psi (\vert \nabla f\vert^2))}\,,
\end{equation}
in which  $f\,:\, \Omega\subset\mathbb{R}^n\to \mathbb{R}$ is any sufficiently smooth function and $\psi$ denotes the inverse(s) of the function $\phi$ defined by \eqref{phi}.
Vector fields  ${\mathbf w}$ expressed by \eqref{sol-k=n-1} or \eqref{altern-k=n-1} are defined on the set $\Omega_f$  except possibly on a singular set $S\subset\gamma_s.$ Often they can be defined with some regularity on $\Omega_f\backslash \gamma_0,$  although  \eqref{vareq-k=n-1} is  not technically satisfied on $\gamma_s$; \emph{cf.} Sec.\ \ref{S-stream2}.

\medskip
In order to obtain solution formulas on non-contractible domains $\tilde \Omega$, one can replace $\nabla f$ in \eqref{sol-k=n-1} and in \eqref{altern-k=n-1} by arbitrary sufficiently smooth vector fields $\mathbf\alpha$ defined on $\tilde\Omega$, satisfying
$\partial_j\alpha_i = \partial_i\alpha_j \;\forall\; i, j.$

\smallskip
The Frobenius condition eq.\ (\ref{frob-k=n-1}) can be obtained for vector functions $\mathbf G$ satisfying the equation
\begin{equation} \label{frobeq-k=n-1}
 \Delta f - \frac{1}{2}\,\nabla f \cdot\,\frac{ \nabla\left(\rho\left(\psi(\vert \nabla f\vert^2)\right)\right)} {\rho\left(\psi\left(\vert \nabla f\vert^2\right)\right)}={\mathbf G}\cdot \nabla f\,,
\end{equation}
or
\begin{equation} \label{frobeqalpha-k=n-1}
 \nabla\cdot \mathbf\alpha - \frac{1}{2}\,\mathbf\alpha \cdot \,\frac{\nabla\left(\rho\left(\psi(\vert \mathbf\alpha\vert^2)\right)\right)}{\rho\left(\psi\left(\vert \mathbf\alpha\vert^2\right)\right)}={\mathbf G}\cdot \mathbf\alpha
\end{equation}
if $\nabla f$ is replaced by $\alpha$ as specified above on a non-contractible domain $\tilde \Omega$.

\medskip
This section extends the results of Secs. \ref{S-stream2}, \ref{S-streamn}, to the system \eqref{vareq-k=n-1}.  The results in the following section, Sec.\ \ref{S-k}, can be used to extend the method to more systems. For an application of the contents of this section, see  Sec.\ \ref{SS-Born-Infeld}.

\section{The general case of  $k$-forms in $n$ dimensions} \label{Hdg}
\label{S-k}
The extension of the method to differential forms of arbitrary order in arbitrary dimension is immediate. In fact, the method can be described naturally in the unifying language of differential forms; see \cite{MO2} for a development in this direction. In this more general context the
nonlinear Hodge--Frobenius equations appear in the form
\begin{equation} \label{A1}
\delta\left(\rho(Q)\omega\right)=0,
\end{equation}
\begin{equation} \label{A2}
d\omega= \Gamma \wedge \omega,
\end{equation}
where we solve for $\omega$, a differential $k$-form on a domain $\Omega\subset\mathbb{R}^n;$ $d:\Lambda^k(\Omega) \rightarrow\Lambda^{k+1}(\Omega)$ is the exterior derivative with formal adjoint $\delta:\Lambda^k(\Omega)\rightarrow\Lambda^{k-1}(\Omega);$ $Q\equiv\vert\omega\vert^2;$  $\rho\,:\, Dom_\rho\subset\mathbb{R}^+\cup\left\{0\right\}\rightarrow\mathbb{R}^+\cup\left\{0\right\}$ is a prescribed, continuously differentiable function on its domain of definition; and
$\Gamma$ is a $1$-form.

Condition \eqref{A2}  guarantees  complete integrability in the cases $k=1$ and $k= n-1$.  If $\Gamma$ can be made exact, say $\Gamma = d\eta$,   $\omega$ is said to be \emph{gradient-recursive} and can be written as
\[
\omega =e^\eta d\zeta
\]
for a $(k-1)$-form $\zeta$; \emph{cf}.\ Sec.\ 2.2 in \cite{MO1}.

If the domain  $\Omega$ is contractible, then applying  the Poincar\'e Lemma as in Secs. \ref{S-stream2}, \ref{S-streamn}, or  \ref{S-k=n-1}, we obtain the solution formula for  $k$-forms $\omega$ satisfying \eqref{A1}  in terms of generalized $\left(n-k-1\right)$-stream forms $f$:
\begin{equation} \label{formsol}
\omega = \frac{\ast df}{\rho\left(\psi\left(\vert df\vert^2\right)\right)}\,,
\end{equation}
where $\ast:\Lambda^k\rightarrow\Lambda^{n-k}$ is the Hodge duality operator on $k$-forms in dimension $n$.
Formula \eqref{formsol} also yields solutions  on non-contractible domains. In that case, one may allow singular stream-forms $f$. Equation \eqref{formsol} may also be replaced by the more general expression
\begin{equation} \label{A4}
\omega = \frac{\ast \alpha}{\rho\left(\psi\left(\vert \alpha\vert^2\right)\right)}
\end{equation}
for  arbitrary closed $(n-k)$-forms $\alpha$.

The Frobenius condition for the $k$ forms $\omega$ expressed by  \eqref{formsol} holds for $1$-forms   $\Gamma$ satisfying
\begin{equation}
\label{frob-A}
\Gamma\wedge\ast df = d\ast df - d\log\rho\left(\psi(\vert df\vert^2)\right)\wedge\ast df\,.
\end{equation}
This is equivalent to the existence of an integrating factor for $\omega$ when $\Gamma$ can be made exact. This can always be done in the case in which $\omega$ is a $1$-form or an $(n-1)$-form.
Equation \eqref{frob-A} can be satisfied with
\begin{equation}\label{Gamma1}
\Gamma = \Gamma_1 - d\log\rho\left(\psi(\vert df\vert^2)\right)\,,
\end{equation}
 whenever $\Gamma_1$ satisfies
\begin{equation}
\label{frob_1-A}
d\ast df = \Gamma_1\wedge\ast df\,.
\end{equation}
This argument shows that the ability to choose $\Gamma$ exact does not depend on the particular type of density function employed, although the value of $\Gamma$ itself does depend on $\rho$. For $\omega$ expressed by \eqref{A4}, formulas \eqref{frob-A}, \eqref{Gamma1} and \eqref{frob_1-A} hold with $df$ replaced by $\alpha$.

\medskip

The equations studied in  Secs.\ \ref{S-stream2} and \ref{S-streamn} correspond to the cases of $1$-forms in dimensions $2$ and $n$, respectively. The system studied in Sec.\ \ref{S-k=n-1} corresponds to the case of $(n-1)$-forms in dimension $n$. In  these cases the differential forms can be easily interpreted as vector functions: a $1$-form can always be interpreted as a vector function, as can an $\left(n-1\right)$-form, being the Hodge-dual of a $1$-form. For the intermediate cases, in which  $k$ is neither $1$ nor $n-1$,  the description in terms of differential forms is more natural. For this reason, an interesting -- perhaps the most interesting ---  illustration of this section is provided by simply taking $k=2,$  $n=4$, and $\rho$ as in \eqref{extr}. This corresponds to the Born--Infeld model in dimension $4$. Our method yields the solutions
\begin{equation}
\label{sol-pm}
\mathbf{\omega_ {\pm}} = \frac {* d f }{\sqrt{\vert d f\vert^2 \pm1}}\,;\;\quad f= \sum_{j=1}^4 a_j(\mathbf x)\, dx^j\in \Lambda^1(\Omega)\,.
\end{equation}
The solutions $\omega_+$ are defined (and uniformly bounded) for smoothly prescribed generalized stream 1-forms $f$, while the solutions  $\omega_-$ require the additional condition $\vert d f\vert>1$ and are unbounded for  choices of  generalized stream forms  which satisfy $\vert d f\vert=1$ at points of the domain $\Omega.$
In this example, the Frobenius condition holds for $1$-forms   $\Gamma$ satisfying
\begin{equation}
\label{frobBI-A}
\Gamma\wedge\ast df = d\ast df - \frac{1}{2} \frac{d\left(\vert df\vert^2\right)}{\vert df\vert^2\pm 1\,}\,.
\end{equation}
See \cite{MO2}, Theorem 3.1  for additional details on the inversion of the function $\phi.$

\section{Examples}
\label{S-exa}
\subsection{A quasilinear elliptic-hyperbolic variational problem}
\label{SS-extremal}
We first illustrate the method on a singular, quasilinear system of variational equations corresponding to (\ref{NLH-intro}, \ref{frob-2}) with $\rho$ given by eq.\ (\ref{extr}). The resulting elliptic-hyperbolic system provides a model for certain nonparametric extremal surfaces, as observed in Sec.\ 1.

Suppose that a surface $\Sigma$ is locally the graph of a function
$\zeta$ with
$\mathbf{w}\equiv\nabla\zeta$ satisfying  eqs.\ (\ref{NLH-intro}, \ref{exact}) with $\rho$ given by (\ref{extr}). Replacing condition (\ref{exact}) with condition (\ref{frob-2}) corresponds to multiplying, at each point of $T\Sigma,$ the length of the gradient vector $\nabla \zeta$ by a
conformal factor $\exp[\eta].$ In this example, the function $\phi$ evaluated at $Q$ is
$$\phi(Q) = \frac{Q}{|1-Q|}\,,\; Q\neq 1$$
and
\begin{equation}\label{dphiextr}
\frac{d\phi(Q)}{dQ}=\pm \frac{1}{(1-Q)^2}\,,
\end{equation}
where the plus sign corresponds to the interval $Q<1$ and the minus sign to the interval  $Q>1,$ respectively.
For $Q<1$, which in the case $n=2,$ $\mathbf{G}=0$ corresponds geometrically to a space-like hypersurface in $\mathbb{M}^3,$ the system of equations (\ref{NLH-intro}, \ref{frob-2}) is elliptic. In that case, the function $\psi_s={[\phi_{|_{[0,1)}}]}^{-1}\,:\, [0,\infty)\to[0,1)$ is given by
\begin{equation}\label{invextr1}
\psi_s:\xi\rightarrow\frac{\xi}{\xi+1}.
\end{equation}
For $Q>1$, which in the case $n=2,$ $\mathbf{G}=0$ corresponds geometrically to a time-like hypersurface in $\mathbb{M}^3,$ the  equations  (\ref{NLH-intro}, \ref{frob-2}) are hyperbolic. In that case the orientation-reversing function $\phi_{|_{(1,\infty)}}$ is again invertible, with inverse given by the orientation-reversing function
$\psi_t=\phi_{|_{(1,\infty)}}^{-1}\,:\, (1,\infty)\to (1,\infty),$
\begin{equation}\label{invextr2}
\psi_t:\xi\rightarrow\frac{\xi}{\xi-1}.
\end{equation}
Specializing to $n=2$ for simplicity, for any given function $f$ one may always use $\psi=\psi_s$ as defined in \eqref{invextr1}, substitute into eq. \eqref{vareqsol}, and obtain a solution to \eqref{NLH-intro}:
\begin{equation} \label{ws}
\mathbf{w}_s = \frac {\nabla_\perp f}{\sqrt{\vert\nabla f\vert^2+1}}.
\end{equation}
(In higher dimensions, $f$ would be replaced by a stream matrix $F,$ as in Sec.\ \ref{S-streamn}.) The solution (\ref{ws}) is at least as regular as $\nabla f$ and remains bounded -- in fact, $Q$ remains in the interval $\left[0,1\right).$ We observe that, although such $\mathbf{w}_s$ is smooth and bounded on $\mathbb R^2,$ $\rho$ becomes singular over the set
\[
\gamma_\infty=\{ (x_0,y_0)\in \mathbb R^2\;:\; \lim_{(x,y)\to(x_0,y_0)}\vert \nabla f\vert^2=\infty\}
\]
corresponding to $Q\to 1.$ Of course this set will be empty provided $f$ is prescribed so that $\nabla f$ is bounded on bounded domains.  If $\gamma_\infty \neq \emptyset,$ the vector function $\mathbf{w}$ may still be defined as a limit. For these solutions $\gamma_0=\emptyset$ and $\gamma_s=\gamma_\infty$; \emph{cf}.\ Sec.\ \ref{S-stream2}.

If $\vert\nabla f(x,y)\vert\in[0,1)$, the use of \eqref{invextr1} to invert $\phi$ is the only possible choice. However, on regions of $\mathbb R^2$ on which $\vert\nabla f(x,y)\vert\in (1,\infty)$ one may use either \eqref{invextr1} or \eqref{invextr2} to obtain solutions to \eqref{NLH-intro}. When the latter choice is made, one obtains
\begin{equation} \label{wt}
\mathbf{w}_t = \frac {\nabla_\perp f}{\sqrt{\vert\nabla f\vert^2-1}}.
\end{equation}
As in the previous case, although such vector functions $\mathbf{w}_t$ remain bounded and may be well-defined on the set $\gamma_\infty$ (because $Q\to 1$ for $\vert \nabla f\vert\to \infty$), eq. \eqref{NLH-intro} is not satisfied at $\gamma_\infty$, as $\rho$ would become singular. In distinction to the previous case, $\mathbf{w}_t$ becomes unbounded if $\vert \nabla f\vert= 1$ on some subset of $\Omega.$ In this case $\Omega_f \subset\Omega$ is a proper inclusion as the set on which $Q=\psi(\vert \nabla f\vert^2)$ becomes unbounded is technically excluded from the definition of  $\Omega_f$; \emph{cf}.\ Sec.\ \ref{S-stream2}.

The choice to invert $\phi$ using \eqref{invextr1} throughout
corresponds to a choice to solve the system (\ref{NLH-intro}, \ref{frob-2})
in the elliptic regime. A solution  is obtained except on the set
$\gamma_\infty.$ Alternatively, one may
choose to use \eqref{invextr2} in one or more regions on which
$\vert\nabla f(x,y)\vert\in (1,\infty).$ It is easy to see that the 1-forms
$\mathbf{w}_s$ and $\mathbf{w}_t$ may only be patched together to make a
continuous vector function $\mathbf{w}$ along curves in $\gamma_\infty,$ as $Q\to 1$ for $(x,y)$ approaching $\gamma_\infty$. The
resulting vector function $\mathbf{w}$ would solve the elliptic-hyperbolic
system (\ref{NLH-intro}, \ref{frob-2}) except at points of $\gamma_\infty.$ In order to avoid singularities for
$\mathbf{w}$, one should choose $f$ so that $\vert\nabla f\vert$ is
bounded away from 1 on bounded regions on which we choose
$\mathbf{w}=\mathbf{w}_t.$

\subsubsection{The patching} \label{secpatch}

As a simple example to illustrate the patching along $\gamma_\infty$ of a solution to system (\ref{NLH-intro}) of type $\mathbf w_s$ with a solution of type $\mathbf w_t$, with $\rho$ given by \eqref{extr}, we consider the radial function
$ f(x,y) = \pm \log \vert r-1\vert.$ Here $r\equiv \sqrt{x^2 + y^2}$,
where the plus sign holds for $r<1$ and the minus sign holds for $r>1$.
This function and its transverse gradient,
$$\nabla_\perp f = \frac{1}{\vert r-1\vert}\left(-\frac{y}{ r} , \frac{x}{ r}\right)\,,$$
are defined except on the set $\{r=1\}\subset\mathbb{R}^2$ and at the origin.
The corresponding solutions $\mathbf w_s$ and $\mathbf w_t$,   expressed by \eqref{ws}, and  \eqref{wt}  respectively,  are
$$
\mathbf w_{s;t} =\frac{1}{\sqrt{1\pm (r-1)^2}}\left(-\frac{y}{ r} , \frac{x}{ r}\right)\,.$$
The solution $\mathbf w_s$ is defined on $\mathbb{R}^2\backslash\{(0,0)\}$, although $f$ was not defined on the set  $\{ r=1\}$, while
$\mathbf w_t$ is defined on $\{(x,y)\in\mathbb{R}^2\,:\, r\in(0,2)\}.$
The two solutions can be considered separately, or one can obtain additional solutions by patching along the circle
$$\gamma_\infty \equiv \{ (x,y)\in \mathbb{R}^2\,:\, r=1\}\,,$$
as on this circle $\mathbf w_s =\mathbf w_t  = \left(-y/{ r} , x/{ r}\right).$ An elementary calculation shows that this can be done with continuity of the first partial derivatives.

For example, one can define $\mathbf w = \mathbf w_t$ on the punctured disk
 \[
 D^\prime\equiv\left\{(x,y)\in \mathbb{R}^2\,:\, r\in \left(0,1\right]\right\},
 \]
 and $\mathbf w  = \mathbf w_s$ on the external region
 \[
 E\equiv\left\{(x,y)\in \mathbb{R}^2\,:\, r\in \left[1,\infty\right)\right\}.
 \]
The resulting vector function
$\mathbf w$ satisfies $\mathbf w \in\mathcal C^1\left( \mathbb{R}^2\backslash\{\left(0,0\right)\}\right).$ Technically
$\mathbf w$ solves  \eqref{NLH-intro}, with $\rho$ given by \eqref{extr} only on the points of its domain for which $r \ne 1.$ A direct calculation employing eqs. \eqref{G-2d} and \eqref{G_1eq-2d} shows that there exist conservative vector fields $\mathbf{G}$ defined on $\mathbb{R}^2\backslash\left\{\left(0,0\right)\right\}$ such that $\mathbf w$ satisfies the Frobenius condition \eqref{frob-2}. One such vector field, corresponding to a radial choice of $\eta$ in (\ref{exactvectorial}), is
\[
\mathbf{G}_{\pm} = \left[\frac{1}{r} - \frac{\left\vert r - 1\right\vert}{1 \pm \left(r-1\right)^2}\right]\left(\frac{x}{r}, \frac{y}{r}\right),
\]
in which the plus sign is to be taken on $E$ and the minus sign is to be taken on $D^\prime.$ Observe that the family $\left\{\mathbf{G}+\lambda \mathbf{w}\right\},$ in which $\lambda$ is a scalar function, also satisfies (\ref{G-2d}, \ref{G_1eq-2d}).  This family contains some additional (nonradial) conservative vector fields. For such fields, the vector $\mathbf w$ constructed above is a solution to the elliptic-hyperbolic system (\ref{NLH-intro}, \ref{frob-2}).

\subsection{Shallow flow with vorticity}
\label{SS-shallow}
Recall from Sec.\ 1 that if we choose $\rho$ as in eq.\ (\ref{shallow}), then we obtain from eqs.\ (\ref{NLH-intro}, \ref{exact}) the continuity equation for shallow, steady, irrotational hydrodynamic flow under a convenient normalization. In this case $\mathbf{w}$ is the flow velocity. The flow potential $\zeta$ exists locally by condition (\ref{exact}) and satisfies
\begin{equation}
\label{stream1}
\rho(Q)\,\partial_x \zeta= -\partial_y f\,,\quad \rho(Q)\, \partial_y \zeta= \partial_x f\,,
\end{equation}
where $f\left(x,y\right)$ is the stream function. That is, in this application the function $f$ of Sec.\ref{S-stream2} has a physical interpretation as an actual stream function; but if (\ref{exact}) is replaced by (\ref{frob-2}), the local flow potential $\zeta$ exists modulo an integrating factor and the methods of, \emph{e.g.}, \cite{Be}, Secs.\ 8 and 9, do not apply in an obvious way. In that case eq. \eqref{NLH-intro} is transformed  locally into the system
$$\rho(Q)\, w_1= -\partial_y f\,,\quad \rho(Q)\, w_2= \partial_x f.$$ Multiplying this equation by $e^{-\eta}$ and imposing the Frobenius condition \eqref{frob-2} with $\mathbf G= \nabla \eta,$
eqs.\ \eqref{stream1} are replaced by
\begin{equation}
\label{stream2}
\rho(Q)\, \partial_x \zeta= -e^{-\eta}\,\partial_y f\,,\quad \rho(Q)\, \partial_y \zeta= e^{-\eta}\,\partial_x f.
\end{equation}
Equations \eqref{stream2}  allow for an interpretation of  $f$ and $\zeta$  as $A$-harmonic functions, that is, as functions satisfying the generalized Cauchy--Riemann equations
\begin{equation}
\label{A-harm-gen}
A(x, \vert \nabla\zeta\vert^2 ) = \nabla_\perp f\,.
\end{equation}
In this case, $A(x, \vert \nabla\zeta\vert^2 ) \equiv e^{\eta(x)}\rho\left(e^{2\eta(x)}\vert\nabla\zeta\vert^2 \right)\nabla\zeta\,.$

\smallskip

Proceeding as in Sec.\ \ref{S-stream2},  we find that the function $\phi$ is strictly increasing on the regime of \emph{tranquil flow}, on which the equations (\ref{NLH-intro}, \ref{frob-2}) are an elliptic system, and strictly decreasing on the regime of \emph{shooting flow}, on which those equations are a hyperbolic system.
More precisely, in this example
the function $\phi(Q)$ is
\begin{equation}\label{phiex}
\phi(Q)= \left(1-\frac{Q}{2}\right)^2 Q
\end{equation}
 and its derivative is
\[
\frac{d\phi(Q)}{dQ} = \left(1-\frac{Q}{2}\right)\left(1-\frac{3}{2} Q\right)\,,
\]
which is strictly positive for $Q\in [0, \frac{2}{3})\cup (2,+\infty)$ and strictly negative for
$Q\in (\frac{2}{3}, 2).$ However, when the squared speed $Q$ exceeds 2, the mass density $\rho$ becomes negative and the physical model no longer applies. Because
\[
\phi \left(\frac{2}{3}\right) = \left(\frac{2}{3}\right)^3
\]
and
\[
\phi(2)=\phi(0) =0,
\]
we have the following correspondences of explicit ranges on which the function $\phi(Q)$ can be inverted:
\begin{align}
\left[0, \frac{2}{3}\right) \to \left[0, \left(\frac{2}{3}\right)^3\right),\label{1}\\
\left(\frac{2}{3}, 2\right) \to \left(0, \left(\frac{2}{3}\right)^3\right).\label{2}\\
\notag\end{align}
If we are willing to consider negative values of $\rho,$ then we also have a third case:
\[
(2, \infty) \to (0, \infty).
\]
The correspondence \eqref{2} is orientation-reversing and corresponds to the hyperbolic regime for the system \eqref{NLH-intro}, \eqref{frob-2}. The explicit expression for the inverse of $\phi$ on each of the three ranges above is given by the cubic formula. Denote by $\psi_1,$ $\psi_2$, $\psi_3$ the three inverses of $\phi$ corresponding to $Q\in [0, \frac{2}{3}),$ $Q\in (\frac{2}{3}, 2)$, and $Q\in (2, \infty),$ respectively. When
\[
\vert\nabla f\vert^2\in \left(0, \left(\frac {2}{3}\right)^3\right),
\]
the function $\psi$ is not uniquely determined, but may be chosen to be $\psi_j$ for  $j=1,2,$ or 3. If
\[
\vert\nabla f\vert^2>  \left(\frac {2}{3}\right)^3,
\]
then the only possible choice is $\psi =\psi_3.$ Obviously, when making these choices care must be taken to obtain reasonably smooth solutions $\mathbf{w}$ which have a plausible interpretation in the context of the chosen model. Note also that the above intervals depend on the constant in Bernoulli's law, which affects the constants in (\ref{shallow}); see the discussion of eq.\ (2.43) in \cite{O4}.

 We will return to this example in Sec.\ \ref{SS-special}.

\subsection{Light near a caustic}
\label{SS-caustic}
If we choose
\begin{equation}\label{optmodel}
\rho(Q)=\sqrt{\left\vert 1 -\frac{\tau^2}{Q}\right\vert},
\end{equation}
then eqs.\ (\ref{NLH-intro}, \ref{exact}) arise in
an elliptic--hyperbolic model for light in the neighborhood of a
caustic with the constant $\tau^2$ representing the index of
refraction; see Sec.\ 1 of \cite{MT}, Sec.\ 3 of \cite{DMV}, or Ch.\
5 of \cite{O4}. The function $\phi$ is invertible on the
\emph{shadow} region of attenuated waves, which corresponds to the
elliptic regime of the system -- that is, to the interval $I_1=
[\tau^2, \infty).$ The function $\phi$ is again invertible on the
\emph{illuminated} region of propagating waves, corresponding
to the hyperbolic regime of the system -- that is, to the interval
$I_2 = [0, \tau^2].$ The boundary between the regions of darkness
and light is the caustic, on which the type of the system
(\ref{NLH-intro}, \ref{exact}) degenerates. As in the preceding cases, replacement of (\ref{exact}) with (\ref{frob-2}) introduces a multiplicative conformal factor in eq.\ (\ref{cons}) which complicates standard techniques, such as quasiconformal mappings, which are based on the existence of a field potential; \emph{cf}.\ \cite{DMV}, Sec.\ 3.

Applying the arguments of Sec.\ \ref{S-stream2}, we obtain the functions
$$\psi_d \,:\, [0,\infty)\to [\tau^2,\infty)\,,\;\xi\to \xi+ \tau^2\,,$$
$$\psi_\ell \,:\, [0,\tau^2]\to [0,\tau^2]\,,\;\xi\to \tau^2-\xi\; \mbox{ (orientation reversing)}\,,$$
yielding
\begin{equation}\label{optics}
\mathbf{w} = \left(\sqrt{\tau^2 \pm\vert\nabla f\vert^2 }\right)\frac{\nabla_\perp f}{\sqrt{\vert\nabla f\vert^2}}\,,
\end{equation}
where the plus sign corresponds to the choice $\psi=\psi_d,$ and the
minus sign, to the choice $\psi=\psi_\ell,$ in \eqref{vareqsol}. We
conclude from this formula that the condition $\vert\nabla
f\vert^2\ne 0$ is not needed for various choices of $f,$ although
relaxing this condition may lead to singularities (in particular, to point singularities when $\vert \nabla f\vert= 0$ at an isolated
point). Again, for any given $f$ one can use $\psi_d$ throughout
and obtain a solution $\mathbf{w}$ (corresponding to the plus sign in
\eqref{optics}), defined everywhere except possibly at points at which
$\vert\nabla f\vert =0;$ one may also decide to use $\psi_\ell$ in
one or more regions of the plane at which $\vert \nabla
f\vert^2\leq\tau^2.$ The two types of solutions can be patched
together only along curves $\gamma$ on which $\vert\nabla f\vert =0$. If $f$ is chosen to satisfy (\ref{frob-2}) for some conservative vector field $\mathbf{G},$ then one obtains in this way a solution to the elliptic--hyperbolic system (\ref{NLH-intro}, \ref{frob-2}).

\subsection{The Born--Infeld model}
\label{SS-Born-Infeld}
In this example we consider the system of equations \eqref{vareq-k=n-1} in $3$ dimensions, with $\rho(Q)$ as in (\ref{extr}).
As computed in Sec.\ \ref{S-k=n-1} for the general case of $n$ dimensions, the solutions ${\mathbf w}$ can be expressed locally by
$${\mathbf w} =\frac{\nabla f }{\rho (\psi (\vert \nabla f\vert^2))}\,,$$
and we will find below explicit forms for the expression $\rho (\psi (\vert \nabla f\vert^2)).$

\smallskip For a physical interpretation of this system with the additional equation \eqref{exact-k=n-1} with $n=3$,
see Sec.\ 1 of \cite{SSY}.  For the effect on the model of using the above equation as a replacement for the more general equation (\ref{frob-k=n-1}) with $n=3,$ see Sec.\ 4 of \cite{MO1}.

\smallskip
This model and the one studied in Section \ref{SS-extremal}, although very different, use the same type of density function $\rho.$ The two models are dimensionally different: the one studied in Section \ref{SS-extremal} is really a model for $1$-forms in $2$ dimensions, while the one studied here is a model for $2$-forms in  $3$ dimensions.  As a consequence, the equations considered are different for the two models. Nonetheless, as the density function $\rho$ and the function $\phi$ are the same in the two models,  much of the analysis carries over to the present model.
So, as already computed in Section \ref{SS-extremal},
$\phi(Q) = Q / |1-Q|$, with $Q\neq 1,$
admitting the inverses $\psi_1\equiv{[\phi_{|_{[0,1)}}]}^{-1}\,:\, [0,\infty)\to[0,1),$ and
$\psi_2\equiv\phi_{|_{(1,\infty)}}^{-1}\,:\, (1,\infty)\to (1,\infty),$ given by
$$
\psi_{1;2}\,:\,\xi\rightarrow\frac{\xi}{\xi\pm1}.$$
We obtain the corresponding solutions
\begin{equation}
\label{w-BI}
\mathbf{w_{1;2}} = \frac {\nabla f }{\sqrt{\vert\nabla f\vert^2 \pm1}}\,.
\end{equation}
The solutions $\mathbf{w_{1}}$ are bounded and  defined for every choice of generalized stream function $f$ (of class $\mathcal C^2$),
while the solutions
$\mathbf{w_2}$ are
defined for generalized stream functions $f$ (of class $\mathcal C^2$) such that $\vert \nabla f\vert^2>1$ and  become unbounded if $\vert\nabla f\vert^2=1$ at points of the domain $\Omega.$
In principle, one could also prescribe a stream function $f$ such that $\vert \nabla f\vert^2\to \infty$ when approaching a surface, say  $\gamma_\infty,$ contained in $\Omega.$  Solutions of type $\mathbf{w_1}$ and of type $\mathbf{w_2}$ could then be glued together along $\gamma_\infty$. Nevertheless, the equations  \eqref{vareq-k=n-1} with $\rho$ prescribed as in \eqref{extr} would no longer be satisfied, as $\rho$ would blow up at the points of $\gamma_\infty.$

\smallskip
In order to find vector functions $\mathbf G$ so that eq.\ (\ref{frob-k=n-1}) will be satisfied as well,  we impose the condition
\[
 \nabla\cdot {\mathbf w} \equiv \frac{\Delta f }{\sqrt{\vert \nabla f\vert^2\pm 1}}- \frac{1}{2}\,\frac{\nabla\left(\vert \nabla f\vert^2\right)\cdot \nabla f }{{\left(\vert \nabla f\vert^2\pm 1\right)}^{\frac{3}{2}}}= \frac{{\mathbf G}\cdot \nabla f}{\sqrt{\vert \nabla f\vert^2\pm 1}}\,.
\]
Multiplying through by $\sqrt{\vert \nabla f\vert^2\pm 1}$, we obtain
an equation for ${\mathbf G}$ of the form
\begin{equation}
\label{G-BI}
\Delta f - \frac{1}{2}\,\frac{ \nabla\left(\vert \nabla f\vert^2\right)}{\vert \nabla f\vert^2\pm 1\,}\,\cdot \nabla f
={\mathbf G}\cdot \nabla f.
\end{equation}
From this we see that  the relation \eqref{frob-k=n-1} can be satisfied with
$$
\mathbf G = - \frac{1}{2}\,\frac{\nabla\left(\vert \nabla f\vert^2\right)}{\sqrt{\vert \nabla f\vert^2\pm 1}}
+ \Delta f\, \frac{\nabla f}{\vert \nabla f\vert^2}\,,
$$
defined on the domain of the solution, with the possible exception of the set
$$
\gamma_G \equiv\{ (x,y,z)\in Dom_{\mathbf w}\subset\mathbb{R}^3\,:\, \nabla f (x,y,z) =\mathbf 0\,; \;\Delta f(x,y,z) \neq 0\,\}\,,$$ in the case of solutions of type $\mathbf w_1$. Note that if $\mathbf G$ is a solution to eq.  \eqref{G-BI} for  given $f$, then the vector functions $\mathbf G + \mathbf H$ with $\mathbf H\cdot\nabla f\equiv 0$ are also solutions to eq. \eqref{G-BI}. One needs to use this degree of freedom to find a conservative $\mathbf G$ satisfying \eqref{frob-k=n-1}.

In particular, choosing $f$ to be the fundamental solution of Laplace's equation in dimension 3
\[
f=f(r)=\frac{1}{r},\,\,r=\sqrt{x^2+y^2+z^2}
\]
yields the solutions
\begin{equation}
\label{fund}
\omega_{1;2}= - \frac { x\, \hat\imath_1 + y\, \hat\imath_2  +z\, \hat\imath_3}{r\sqrt{1 \pm r^4}}\,,
\end{equation}
which satisfy the Frobenius condition \eqref{frob-k=n-1} with
$$\mathbf G = 2 \,\frac { x\, \hat\imath_1 + y\, \hat\imath_2  +z\, \hat\imath_3}{r^4\sqrt{1 \pm r^4}}
+ \mathbf H\,.$$ The vector field $\mathbf G$ can be made conservative by a particular choice of $\mathbf H$ satisfying $\mathbf H\cdot\nabla f\equiv 0.$

As observed earlier for the general solutions $\mathbf w_2$, the solution \eqref{fund} corresponding to the minus sign becomes unbounded on the set $$\{(x,y,z)\in\mathbb{R}^3\,:\, \vert\nabla f\vert^2(x,y,z)\equiv\frac{1}{r^4} =1\}\,,$$ that is, on the sphere $r=1.$
This can often happen when inverting in the hyperbolic regime. In fact, in that case the map $\phi$ is orientation-reversing and $Q$ can become unbounded at points on which the prescribed $\vert \nabla f\vert$ is bounded, and thus on a bounded domain.

\subsection{Special classes of streamlines}
\label{SS-special}
By prescribing certain symmetries in the stream function $f\left(x,y\right),$ one obtains explicit examples of solutions $\mathbf{w}$ to \eqref{NLH-intro}.
For example, assume that the function $f$ can be expressed as
$f(x,y) = \tilde f (t)$, where $t= g\left(x,y\right)$ is a given continuously differentiable function of the given domain $\Omega.$
We obtain
\[\partial_x f = \tilde{f}^\prime(t) \,\partial_x g\,,\qquad \partial_y f = \tilde{f}^\prime(t) \,\partial_y g\,,
\] and
\[
\vert\nabla f\vert^2 = \left[\tilde{f}^\prime(t)\right]^2\vert \nabla g\vert^2.
\]
For any choice of $\rho$, eq. \eqref{normw} yields a solution $\mathbf{w}$ of \eqref{NLH-intro} satisfying
\[
\frac{\mathbf{w}}{\vert \mathbf{w}\vert}=\pm\left(\frac{\nabla_\perp g}{\vert\nabla g\vert}\right).
\]
A special case is the example of radial $f$ -- that is, $t= g\left(x,y\right)\equiv x^2+y^2.$ Then for any choice of $\rho$ we obtain the simple formula
\begin{equation} \label{frad}
\frac{\mathbf{w}}{\vert \mathbf{w}\vert}=\pm \frac{-y\,\hat{\imath}+x\,\hat{\jmath}}{\sqrt t}.
\end{equation}
In this case formula (\ref{gradf}) becomes
\[
\phi(Q)=\vert\nabla f\vert^2 = 4t\tilde f^\prime\left(t\right)^2
\]
and \eqref{altern} yields
\[ \mathbf{w}=\pm \frac{-y\,\hat{\imath}+x\,\hat{\jmath}}{\sqrt t}\, \sqrt{\psi(4t\tilde f^\prime\left(t\right)^2)}.
\]

In applications to the shallow-water model corresponding to
(\ref{shallow}), one may choose the function $f$ so that
$\vert\nabla f\vert^2\left(t\right)$ is the function $\phi(t)$
modulo a constant dilation of the domain. In that case,
$\psi\left(\vert\nabla f\vert^2\right)$ is essentially the identity
on the domain. We illustrate this application with a numerical example. Choosing
\[
f\left(x,y\right) = \left(t -
\frac{t^2}{4R}\right)\frac{1}{2\sqrt{R}},
\]
where $R$ is a given positive constant, we have
\[
\vert\nabla
f\vert^2\left(t\right)=\left(1-\frac{t}{2R}\right)^2\frac{t}{R}.
\]
We obtain, as expected, a solution having the familiar form
\[
\mathbf{w} =
\frac{\frac{-y}{\sqrt{R}}\left(1-\frac{t}{2R}\right)\,\hat{\imath}+\frac{x}{\sqrt{R}}\left(1-\frac{t}{2R}\right)\,\hat{\jmath}}{1-\psi\left(\frac{t}{R}\left(1-\frac{t}{2R}\right)^2\right)/2}=
\]
\[
\frac{\frac{-y}{\sqrt{R}}\left(1-\frac{t}{2R}\right)\,\hat{\imath}+\frac{x}{\sqrt{R}}\left(1-\frac{t}{2R}\right)\,\hat{\jmath}}{1-\frac{t}{2R}}=\frac{-y\,\hat{\imath}+x\,\hat{\jmath}}{\sqrt{R}},
\]
where we have applied the various inverse functions $\psi$ described
in Sec.\ref{SS-shallow}  on the various subdomains $B_{2R/3},$ $B_{2r}\backslash
B_{2R/3},$ and $\mathbb{R}^2\backslash B_{2R}.$ In the second to last equality
we have used the identity
\[
\psi\left(\rho^2\left(\frac{t}{R}\right)\frac{t}{R}\right)=\psi\circ\phi\left(\frac{t}{R}\right)=\frac{t}{R}\,,\qquad t=x^2+y^2\,,
\]
which is satisfied by all the inverses.

The vector $\mathbf{w}$ so defined satisfies (\ref{NLH-intro}) with $\rho$
given by (\ref{shallow}) -- that is, with
\[
\rho = 1 -\frac{t}{2R},
\]
and (\ref{frob-2}) with
\[
\mathbf{G} = \frac{2x}{t\sqrt R}\,\hat{\imath}+\frac{2y}{t\sqrt R}\,\hat{\jmath}
\]
on $\mathbb{R}^2\backslash \left\{\left(0,0\right)\right\},$ and is
globally smooth. The system (\ref{NLH-intro}, \ref{frob-2}) is elliptic on $B_{2R/3}\backslash \{(0,0)\},$
hyperbolic on $B_{2r}\backslash
B_{2R/3},$ and then again elliptic on $\mathbb{R}^2\backslash B_{2R}$ (although that region would correspond to negative $\rho$).

More examples can be found by applying the same procedure to a more general function $t=g(x,y)$, such that $\vert \nabla g\vert$ would be a function only of $t$.

\medskip

\noindent \textbf{Acknowledgment}. The authors are grateful to Yuxi Zheng for discussion of the gas-dynamics research literature.


\begin{thebibliography} {99}

\bibitem{Be} L. Bers, \emph{Mathematical Aspects of Subsonic and
Transonic Gas Dynamics}, Wiley, New York, 1958.

\bibitem{Ch} C. J. Chapman, \emph{High Speed Flow}, Cambridge University Press, Cambridge, 2000.

\bibitem{CF} R. Courant and K. O. Friedrichs, \emph{Supersonic Flow and Shock Waves}, Springer-Verlag, Berlin-Heidelberg-New York, 1976.

\bibitem{DMV} E. De Micheli and G. A. Vianoz, The evanescent waves in geometrical optics and
the mixed hyperbolic-elliptic type systems, \emph{Applicable Analysis} \textbf{85}
(2006), 181-204.

\bibitem{E}  D. G. B. Edelen, \emph{Applied Exterior Calculus}, Wiley, New York, 1985.

\bibitem{ISS} T. Iwaniec, C. Scott, and B. Stroffolini, Nonlinear Hodge
theory on manifolds with boundary, \emph{Annali Mat. Pura Appl.}
\textbf{177} (1999), 37-115.

\bibitem{MT} R. Magnanini and G. Talenti, On complex-valued
solutions to a 2D eikonal equation. Part one: qualitative
properties, \emph{Contemporary Math.} \textbf{283} (1999), 203-229.

\bibitem{MO1} A. Marini and T. H. Otway, Nonlinear Hodge--Frobenius equations and the Hodge--B\"acklund transformation, \emph{Proc. R. Soc. Edinburgh} \textbf{140A} (2010), 787-819.

\bibitem{MO2} A. Marini and T. H. Otway, Duality methods for a class of quasilinear systems, \emph{Ann. Inst. Henri Poincar\'{e} C. Nonlinear Analysis.} (in press).

\bibitem{Me} R. E. Meyer, \emph{Introduction to Mathematical Fluid Dynamics}, Dover, New York, 1982.

\bibitem{O1} T. H. Otway, Nonlinear Hodge maps, \emph{J. Math. Phys.} \textbf{41} (2000), 5745-5766.

\bibitem{O2} T. H. Otway, Maps and fields with compressible density, \emph{Rend. Sem. Mat. Univ. Padova} \textbf{111} (2004), 133-159.

\bibitem{O4} T. H. Otway, \emph{The Dirichlet Problem for Elliptic-Hyperbolic Equations of Keldysh Type}, Lecture Notes in Mathematics, Vol. 2043, Springer-Verlag, Berlin-Heidelberg-New York-Tokyo, 2012.

\bibitem{SS1} L. M. Sibner, R. J. Sibner, A nonlinear Hodge-de Rham theorem, \emph{Acta Math}. \textbf{125} (1970), 57-73.

\bibitem{SSY}  L. M. Sibner, R. J. Sibner, and Y. Yang, Generalized Bernstein property and gravitational strings in Born-Infeld theory, \emph{Nonlinearity} \textbf{20} (2007), 1193-1213.

\bibitem{St} J. J. Stoker, \textit{Water Waves}, Interscience, New York, 1987.

\end{thebibliography}
\end{document}